%% file: root.tex
\documentclass[10pt]{article}

\def\figpath{figures}
\input{\texpath/preamble}

\begin{document}
\maketitle
\forECC{
	\thispagestyle{empty} 
	\pagestyle{empty}
}
\begin{abstract}                
	This article concerns the performance limits of strictly causal state estimation for linear systems with fixed, but uncertain, parameters belonging to a finite set. In particular, we provide upper and lower bounds on the smallest achievable gain from disturbances to the point-wise estimation error. The bounds rely on forward and backward Riccati recursions---one forward recursion for each feasible model and one backward recursion for each pair of feasible models. We give simple examples where the lower and upper bounds are tight. 
\end{abstract}

\input{\texpath/introduction}
\input{\texpath/notation}
\input{\texpath/body}
\section{ACKNOWLEDGEMENTS}
The author expresses sincere gratitude to colleagues Anders Rantzer and Venkatraman Renganathan for their stimulating and insightful discussions on the content and organization of this paper.
\input{\texpath/appendix}

\forECC{
	\bibliographystyle{IEEEtran}
}
\forArxiv{
	\bibliographystyle{plainnat}

}
\bibliography{references}             

\end{document}

%% file: tex/preamble.tex
\newcommand{\forArxiv}[1]{#1}  
\newcommand{\forECC}[1]{}      

\usepackage{verbatim, color}
\usepackage{amsmath,amssymb, amsfonts}

\usepackage{amsthm}
\usepackage{url}

\usepackage{booktabs}
\usepackage{graphicx}
\usepackage{balance}
\usepackage{tikz, pgfplots}
\usepackage{subcaption}
\usepackage{mathtools} 
\forArxiv{
	\usepackage{natbib}
	\renewcommand{\cite}{\citep}
	\usepackage[titletoc,title]{appendix}
}

\DeclareMathOperator{\bdiag}{BDiag}

\newcommand{\onote}[1]{{\color{black}#1}}
\newcommand{\bmat}[1]{\begin{bmatrix} #1 \end{bmatrix}}

\newcommand{\R}{\mathbb{R}}
\newcommand{\modelset}{\mathcal{M}}
\newcommand{\seq}[3]{#1_{[#2:#3]}}
\def \nmodels {M}

\newcommand*{\defeq}{\stackrel{\small{\mathsf{def}}}{=}} 
\newcommand{\tran}{{\mkern-1.5mu\mathsf{T}}}
\newcommand{\stack}[1]{\mathbf{#1}}

\newtheorem{theorem}{Theorem}
\newtheorem{prop}{Proposition}
\newtheorem{lemma}{Lemma}

\newtheorem{remark}{Remark}

\mathtoolsset{centercolon}

\usetikzlibrary{positioning, shapes, arrows, calc }
\tikzstyle{block} = [draw, rectangle, line width=0.7mm]
\tikzset{sumcircle/.style={draw,circle,line width=1.4pt,outer sep=0pt,label=center:{$+$},minimum width=2em}}
\tikzset{prodcircle/.style={draw,circle,line width=1.4pt,outer sep=0pt,label=center:{$\times$},minimum width=1em}}

\forECC{
	\title{\LARGE \bf Minimax Performance Limits for Multiple-Model Estimation
	\thanks{*This project has received funding from the European Research Council (ERC) under the European Union's Horizon 2020 research and innovation programme under grand agreement No 83142 (ScalableControl).}%
	\author{\authorblockN{Olle Kjellqvist}
		\authorblockA{Department of Automatic Control, Lund University, Sweden\\
			Email: {\tt\small olle.kjellqvist@control.lth.se}
		}
	}
}
}
\forArxiv{
\title{ Minimax Performance Limits for Multiple-Model Estimation
	\thanks{This project has received funding from the European Research Council (ERC) under the European Union's Horizon 2020 research and innovation programme under grand agreement No 83142 (ScalableControl).}%
	\author{Olle Kjellqvist
		\thanks{Department of Automatic Control, Lund University, Sweden. 
			Email: {\tt\small olle.kjellqvist@control.lth.se}
		}
	}
} 
}

%% file: tex/introduction.tex
\section{Introduction}

Multiple-model estimation is a valuable tool for state estimation of systems that operate in different modes, for problems involving unknown parameters, for dealing with systems subject to faults, and for target tracking.
If the mode is known, one selects the filter corresponding to the current mode.
Otherwise, one can use a bank of filters, one for each mode, and cleverly combine the estimates.
The latter approach is precisely what is called multiple-model estimation.

Almost all of the literature assumes that the system is affected by stochastic noise and that good noise statistics are available.
Unfortunately, many popular methods are sensitive to a mismatch between the assumed and actual noise statistics.
This assumption limits the applicability of in control systems, where we often use simplified models and disguise the model mismatch as additive disturbances.
These disturbances are sometimes poorly modeled by Gaussian noise, and the noise statistics are often unknown.

\begin{figure*}
	\centering
	\forECC{\includegraphics[width = .8\textwidth]{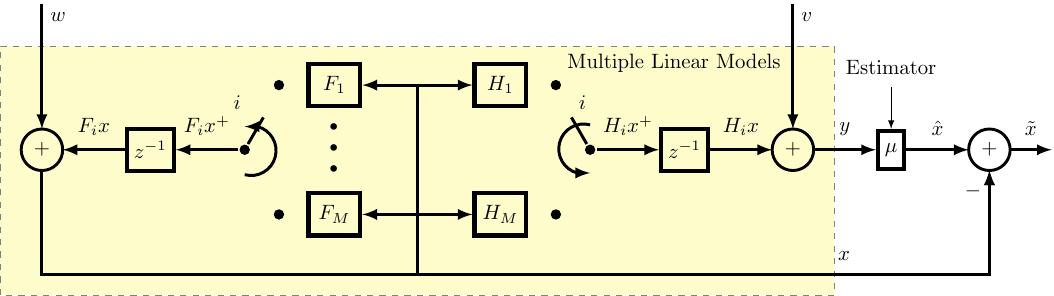}}
	\forArxiv{\includegraphics[width = \textwidth]{figures/newprob.pdf}}
	\caption{An illustration of the multiple-model estimation problem.}\label{fig:problem}
\end{figure*}

In this article, we consider the problem of predicting the state of a linear system with unknown but fixed parameters belonging to a finite set.
We assume that the system is affected by disturbances but make no assumptions about the noise statistics.
We study the \emph{minimax performance level}, defined as the gain from disturbances to point-wise estimation error, and are concerned with bounding the optimal (smallest achievable) performance level.
See Fig.~\ref{fig:problem} for an illustration of our problem.

\subsection{Contributions}
This author, and Rantzer, recently proposed an estimator that achieves the optimal performance level but the performance level itself was not characterized \cite{Kjellqvist2022Minimax}.
The main contribution of this article is to extend the framework in \cite{Kjellqvist2022Minimax} with a method to compute upper and lower bounds of the optimal performance level.
These bounds are computed offline, a priori, and depend on the pairwise interaction between candidate models.

\subsection{Background}
The idea of using multiple models to reduce uncertainty is prevalent in many fields.
It has been used in adaptive estimation since the '60s~\cite{Magill1965}, where it is called multiple-model estimation and in feedback control since the '70s~\cite{Athans1977MMAC}, where it is called multiple-model adaptive control~\cite{Buchstaller2016Framework}, or supervisory control~\cite{Hespanha2001Tutorial}. 
The concept has been known in machine learning at least since Dasarty and Sheela introduced the ``Composite classifier system'' in 1979~\cite{Dasarty1979Composite}, and is commonly referred to as ensemble learning~\cite{Dietterich2000EnsembleMI}.
In the field of economics, the idea of multiple models is known as model averaging~\cite{Steel2020ModelAveraging}, and was popularized by the work of Bates and Granger~\cite{Bates1969TheCO}. 

The task usually falls into one of two categories: \emph{model selection}, where the goal is to find the best performing model, or \emph{model averaging}, where the goal is to use all the models to generate an estimate of some common quantity.
In this article, the focus is on predicting the state in dynamical systems, which falls into the latter category.

When the model is known, the Kalman filter is the realization of many reasonable estimation strategies.
The minimum variance estimate, the maximum-likelihood estimate, and the conditional expectation under white-noise assumptions~\cite{Anderson79} all coincide with the estimate generated by the Kalman filter.
The filter also has appealing deterministic interpretations as the minimum energy estimate~\cite{Willems2004Deterministic, Buchstaller2020Deterministic}, and as Krener showed~\cite{Krener1980Minimax}, it constitutes a minimax optimal estimate.

Interestingly, a minimax optimal estimate can be derived and computed without explicit knowledge of its \emph{minimax performance level}, a property not shared with the $\mathcal H_\infty$-optimal estimate~\cite{Shen97} and controller~\cite{Basar95}, which require knowledge of their performance levels. Tamer Ba\c{s}ar showed that the optimal performance level can be obtained from the finite escape times of some related Riccati recursions~\cite{Basar1991Optimum}.

In the case of multiple fixed models, the different estimation strategies give rise to different estimates\footnote{Except the maximum likelihood estimate under white-noise assumptions and a uniform prior over $\mathcal M$ coinciding with the minimum-energy estimate. }. 
The stochastic multiple-model approach to adaptive estimation was introduced in the '60s~\cite{Magill1965, Lainiotis1971} for linear systems with fixed, but unknown parameters, and has numerous applications in fault detection, state estimation and target tracking ~\cite{Rong2005Survey}.
This estimation algorithm applies the Bayes rule recursively under white-noise assumptions on $(w,v)$ and is well described in many textbooks like \cite{Gustafsson2000, Crassidis2011, Anderson79}. 
The book \cite{Anderson79} also contains a convergence result, stating that given a certain distinguishability condition\footnote{Silvestre et al., \cite{Silvestre2020Distinguishability}, recently reexamined the distinguishability requirements from a multiple-model adaptive control perspective.}, the conditional probability for the active model generating the data converges to $1$ as time goes to infinity. 
Vahid et al., \cite{Vahid2009Energy}, proposed a minimum-energy condition for multiple-model estimation and proved a convergence result given a persistency-of-excitation-like criterion.

Multiple-model estimation has also been extended to the case with changing parameters, the case when $i$ in Fig.~\ref{fig:problem} evolves on a Markov chain. One can, in principle, solve exactly for the Baysian average, but this is computationally intractable as the number of feasible trajectories grow exponentially with time.
Instead, there exist sub-optimal algorithms that cleverly combine estimates at each time-step, compressing the feasible trajectories, like Blom's Interacting-Multiple-Model algorithm, \cite{Blom1988}.
This idea was further generalized by Li and Bar-Shalom to the case when the model set varies with time, \cite{Li1996Variable}.

The work in this article is inspired by recent progress in \emph{minimax adaptive control}~\cite{Rantzer2020, Rantzer2021, Kjellqvist2022Learning}, and in a broader sense, the search for performance guarantees in learning-based control and identification~\cite{Matni2019Reinforcement, Mania2022TV}.

\subsection{Outline}
The rest of this paper is organized as follows.
We establish notation in Section~\ref{sec:notation}.
Section~\ref{sec:minimax} contains the problem formulation and solution.
Illustrative examples are in Section~\ref{sec:examples}.
We give conclusions and final remarks in Section~\ref{sec:conclusions}.
The proofs of the main results and supporting Lemmata are contained in the appendix.

%% file: tex/notation.tex
\section{Notation}
\label{sec:notation}
The set of $n \times m$-dimensional matrices with real coefficients is denoted $\R^{n\times m}$. 
The transpose of a matrix $A$ is denoted $A^\tran$. 
For a symmetric matrix $A \in \R^{n\times n}$, we write $A \succ (\succeq) 0$ to say that $A$ is positive (semi)definite. 
The $n\times n$-dimensional identity matrix is denoted $I_n$, and the $n \times m$-dimensional zero matrix is denoted $0_{n\times m}$.
Given $x\in \R^n$ and $A\in \R^{n\times n}$, $|x|^2_A := x^\tran A x$. 
For a vector $x_t \in \R^n$ we denote the sequence of such vectors up to time $t$ by $x_{[0:t]} := (x_k)_{k=0}^{t}$.
For a sequence of square matrices $(A_i)_{i=1}^M$, we denote the corresponding block-diagonal matrix as $\bdiag(A_i)_{i=1}^M$.

%% file: tex/body.tex
\section{Minimax performance limits}\label{sec:minimax}
\subsection{Problem statement}
In this article, we consider strictly causal\onote{\footnote{\onote{The ideas in this paper extend to other information structures like filtering, $k$-step prediction, and smoothing, but they require some extra steps.}}} state estimation for uncertain linear systems of the form
\begin{equation}
        \begin{aligned}
        \label{eq:dynamics}
		x_{t+1} & = Fx_t + w_t, & (F, H) \in \mathcal M \\
		y_t & = Hx_t + v_t,&  0 \leq t \leq N-1, \\
        \end{aligned}
\end{equation}
where $x_t\in \R^n$, and $y_t\in \R^m$ are the states and the measured output at time $t$.
$w_t\in\R^n$ and $v_t\in\R^m$ are unmeasured process disturbance and measurement noise.
We employ a deterministic framework and make no assumptions on the distributions of $w_t$ and $v_t$.
Instead, they are adversarially chosen to maximize the objective of a related minimax problem that we will define shortly.
The model, $(F, H) \in \R^{n \times n} \times \R^{m \times n}$ is unknown but fixed, belonging to a (known) finite set 
\begin{equation*}\label{eq:model_set}
	\mathcal M = \{(F_1, H_1), \ldots ,(F_M, H_M)\}.
\end{equation*}

The state estimate at time $N$, $\hat x_N$, is generated by a causal estimator, $\mu$, that depends on previous measurements but is unaware of the model, $(F, H)$, and noise, $(w,v)$, realizations,
\begin{equation*} \label{eq:mu}
    \hat x_N = \mu(y_{N-1}, \ldots, y_0).
\end{equation*}

We are interested in describing the smallest $\gamma_N$, denoted $\gamma_N^\star$, such that the below expression has finite value.
\begin{multline}\label{eq:def_J}
    J_N^\star(\hat x_0) := \inf_\mu\sup_{x_0, \seq{w}{0}{N-1}, \seq{v}{0}{N-1}, i} \Bigg\{
        |x_N - \hat x_N|^2 \\
         - \gamma_N^2\bigg(|x_0 - \hat x_{0, i}|^2_{P_{0, i}} + \sum_{t=0}^{N-1}\Big[|w_t|^2_{Q_i^{-1}} + |v_t|^2_{R_i^{-1}} \Big] \bigg)
        \Bigg\},
\end{multline}
where the trajectory $\seq{x}{0}{N}$ in \eqref{eq:def_J} is generated according to \eqref{eq:dynamics} with $(F_i, H_i) \in \mathcal M$.
The problem set-up is a two-player game where the adversary picks the disturbance sequences $\seq{w}{0}{N-1}$ and $\seq{v}{0}{N-1}$, the initial state $x_0$, and the active model $i = 1, \ldots, M$. 
The minimizing player picks the estimation policy $\mu$.
The matrices $Q_i \in \R^{n \times n}$ and $R_i \in \R^{m \times m}$ are positive definite matrices that weights the norms on $w$ and $v$.
The matrices $P_{0, i} \in \R^{n\times n} $ are positive definite and quantify the uncertainty in the estimates of the initial states $\hat x_{0, i}$.
\subsection{Forward Recursions}
The forward recursions describe the worst-case disturbances consistent with the dynamics and an observed trajectory. They are also fundamental in constructing a minimax-optimal estimator $\mu_\star$.
The recursions are equivalent to those of a Kalman filter of a system driven by zero-mean independent white noise sequences $w_t$ and $v_t$ with covariance matrices $Q_i$ and $R_i$ respectively,
\begin{equation}\label{eq:Kalman_gain}
    \begin{aligned}
K_{t,i}     & = F_iP_{t,i}H_i^\tran(R_i + H_i P_{t,i} H_i^\tran)^{-1}, \\
P_{t+1,i} & = Q_i + F_iP_{t, i}F_i \onote{-} K_{t, i}(R_i + H_i P_{t,i} H_i^\tran) K_{t, i}^\tran.
    \end{aligned}
\end{equation}
The relation between the stochastic interpretation and our deterministic framework lies in that the least-squares estimate coincedes with the maximum-likelihood estimate under white-noise assumptions. 

\begin{remark}
$P_{0, i}$ is a regularization term that penalizes deviations from an initial state estimate $\hat x_{0,i}$ and can be interpreted as the covariance of the initial estimate $\hat x_{0, i}$. 
It is practical to choose $P_{0, i}$ as the stationary solution to \eqref{eq:Kalman_gain}, and we will do so in the sequel to simplify the notation by removing the time index.
The results in this section are valid for any positive semi-definite choice of $P_{0, i}$.
However, the resulting observer dynamics will be time-varying.
We leave it to the reader to reintroduce the dependence on $t$.
\end{remark}

The solution, $P_i$, to the Riccati equation \eqref{eq:Kalman_gain} quantifies the uncertainty of the state estimate given the observations $y_{0:t}$ and the model $i$ and bounds the smallest achievable gain from below if the model is known. 
This is formalized in the following proposition, whose proof is in the Appendix.
\begin{prop}\label{prop:bound_P}
    $\gamma_N \geq \gamma_N^\star$ only if $P_i \preceq \gamma_N^2 I$ for all $i = 1, \ldots, M$.
\end{prop}

In our previous work, \cite{Kjellqvist2022Minimax}, we show how to construct the minimizing argument $\mu_\star$ of \eqref{eq:def_J} in the case of output-prediction. 
The estimator uses the forward recursions \eqref{eq:Kalman_gain} and requires a $\gamma_N$ that fulfills Proposition~\ref{prop:bound_P}.
The following proposition shows how to construct a state predictor that is optimal for~\eqref{eq:def_J}.
\begin{prop}[Minimax multiple-model estimator] \label{prop:minimax_optimal}
	Given matrices $F_i \in \R^{n \times n}$ and $H_i \in \R^{m \times n}$, positive definite $Q_i,\ P_{0, i} \in \R^{n \times n}$ and $R_i \in \R^{m \times m}$ for $i = 1, \ldots, M$. 
	With $P_{0, i}$, $P_{i}$ and $K_{i}$ as the stationary solutions to \eqref{eq:Kalman_gain},
\begin{equation}~\label{eq:Rtilde}
    \tilde R_i = R_i + H_i P_iH_i^\tran,
\end{equation}
a quantity $\gamma_N$ such that $\gamma_N^2 I \succ P_{i}$, the below estimate achieves the infimum in~\eqref{eq:def_J}:
	\[
		\hat x_N^\star = \min_{\hat x_N}\max_i \left\{|\hat x_N - \breve x_{N, i}|^2_{(I - \gamma^{-2}P_{i})^{-1}} - \gamma_N^2 c_{N, i} \right\},
	\]
	where $\breve x_{N,i} \in \R^n$ and $c_{N,i} \in \R$ are generated according to
	\begin{subequations} \label{eq:estimators}
		\begin{align}
			\breve x_{0,i} & = x_0, \quad c_{0,i} = 0,  \\
			\breve{x}_{t+1,i}   & = F_i \breve{x}_{t,i} + K_{i}(y_t - H_i\breve{x}_{t,i}), \label{eq:estimators:state_update}\\
			c_{t+1,i} & = |H_i\breve x_{t,i} -y_t|^2_{\tilde R_{i}^{-1}} + c_{t,i}.
		\end{align}
    \end{subequations}
\end{prop}
\begin{proof}
	The proof is identical to that of Theorem 1 in~\cite{Kjellqvist2022Minimax} but with the following modifications: $P_{0, i}$ is replaced by the stationary solution to \eqref{eq:Kalman_gain} leading to $K_{t, i}$ and $P_{t, i}$ being replaced by $K_i$ and $P_i$, the term $\hat y_N - H_ix_N$ is replaced by $\hat x_N - x_N$.
\end{proof}
\subsection{Backward Recursions}
The backward recursions are similar to those of the linear-quadratic regulator and relate to the worst-case trajectories, in contrast to the forward recursions, which relate to the worst-case disturbances consistent with any given trajectory. 
They play no role in constructing the optimal estimator, $\mu_\star$, once a performance level $\gamma$ has been found, but form the basis for \emph{a priori} analysis of the optimal performance level $\gamma_N^\star$ that holds for any realization.
Let
\begin{equation*}\label{eq:FK}
	F^{ij} = \bmat{F_i - K_iH_i & 0_{n \times n} \\ 0_{n \times n} & F_j - K_jH_j}, \quad
    K^{ij} = \bmat{K_i \\ K_j}. \quad 
\end{equation*}
$F^{ij}_t$ corresponds to the closed-loop of a pair $(i, j)$ of Kalman filters with filter gains $K_i$ and $K_j$ as in \eqref{eq:Kalman_gain}.
We will express the necessary and sufficient conditions using the following Riccati recursions.
Given some symmetric matrix $T^{ij}_N \in \R^{2n \times 2n}$ and $t = N-1, \ldots, 0$,
\begin{gather} \label{eq:riccati}
\begin{aligned}
	X^{ij}_{t} & = (K^{ij})^\tran T^{ij}_{t + 1}K^{ij} +(\tilde R_i^{-1} + \tilde R^{-1}_j),  \\
		   L^{ij}_{t} & =(X^{ij}_{t})^{-1}\left((K^{ij})^\tran T^{ij}_{t + 1}F^{ij} - \bmat{\tilde R_i^{-1}H_i & \tilde R_j^{-1}H_j} \right), \\
		   T^{ij}_{t} & = (F^{ij})^\tran T^{ij}_{t+1}F^{ij} - (L^{ij}_{t})^\tran X^{ij}_{t}L^{ij}_{t} \forECC{\\ }
		   \forECC{& \qquad} +\bmat{H_i^\tran\tilde R_i^{-1}H_i \\ & H_j^\tran\tilde R_j^{-1}H_j}. 
\end{aligned}
\raisetag{2.5\baselineskip}
\end{gather}
For these recursions to be well-defined, the matrix $X_t^{ij}$ must be invertible.
The conditions for bounding $\gamma_N^\star$ are related to the positive definiteness of $X_t^{ij}$ and are summarized in Theorems~\ref{thm:suff} and~\ref{thm:nec} below. 
The first concerns sufficient conditions and can be used to obtain upper bounds. 
\begin{theorem}[Sufficient Condition] \label{thm:suff}
	Given matrices $F_i \in \R^{n \times n}$ and $H_i \in \R^{m \times n}$, positive definite $Q_i \in \R^{n \times n}$ and $R_i \in \R^{m \times m}$ for $i = 1, \ldots, M$. 
	Further, let $P_{0, i} = P_i$ and $K_i$ be the stationary solutions to \eqref{eq:Kalman_gain}, and consider a quantity $\gamma_N$ such that $\gamma_N^2 I \succ P_i$. 
	Let $\underline Q \in \R^{n \times n}$ be a positive definite matrix such that $\underline Q \preceq I - \gamma_N^{-2} P_i$ for all $i = 1, \ldots, M$ and initialize the backward recursions~\eqref{eq:riccati} with the terminal state
    \[
	    T^{ij}_N = -\bmat{\underline Q^{-1} & -\underline Q^{-1} \\ -\underline Q^{-1} & \underline Q^{-1}}/\gamma_N^2.
    \]    
    Assume that $X^{ij}_t$ in \eqref{eq:riccati} is negative definite for all $i, j$. 
    Then $\gamma_N^\star \leq \gamma_N$ and 
    \[
	    J_N^\star(\hat x_0) \leq \frac{1}{2} \max_{i, j}\left\{-\gamma_N^2\bmat{
            \hat x_{0, i} \\ 
            \hat x_{0, j}
            }^\tran T^{ij}_0 \bmat{
                \hat x_{0, i} \\ 
                \hat x_{0, j}
}\right\}.
    \]
\end{theorem}

The second theorem concerns necessary conditions and helps obtaining lower bounds.
\begin{theorem}[Necessary Condition]\label{thm:nec}
	Given matrices $F_i \in \R^{n \times n}$ and $H_i \in \R^{m \times n}$, positive definite $Q_i \in \R^{n \times n}$ and $R_i \in \R^{m \times m}$ for $i = 1, \ldots, M$. 
	Further, let $P_{0, i} = P_i$ and $K_i$ be the stationary solutions to \eqref{eq:Kalman_gain}, and consider a quantity $\gamma_N$ such that $\gamma_N^2 I \succ P_i$. 
	Initialize the backward recursion \eqref{eq:riccati} with the terminal state
    \[
	    \begin{aligned}
		    T^{ij}_N & = -\bmat{Q^{ij} & -Q^{ij} \\
	-Q^{ij} &Q^{ij}
	}/\gamma_N^2, \\
	Q^{ij} & = (2I - \gamma_N^{-2}(P_{N, i} + P_{N, j}))^{-1}.
\end{aligned}
    \]
    If $X^{ij}_t \not \preceq 0$ for some pair $i, j$ and $0 \leq t \leq N-1$, then $\gamma_N^\star > \gamma_N$.
    If $X^{ij}_t \succ 0$, for all $t = 0, \ldots, N-1$ then
    \[ 
	    J_N^\star(\hat x_0) \geq \frac{1}{2}\max_{ij}\left\{-\gamma_N^2\bmat{\hat x_{0, i} \\ \hat x_{0, j}}^\tran T^{ij}_0 \bmat{\hat x_{0, i} \\ \hat x_{0, j}}\right\}.
    \]
\end{theorem}
\begin{remark}
	Theorems~\ref{thm:suff} and \ref{thm:nec} give upper and lower bounds on $J_n^\star$ that can be translated upper and lower bounds on $\gamma_N^\star$ by bisecting over $\gamma_N$.
\end{remark}
\section{Examples}\label{sec:examples}
Figures~\ref{fig:unstablepm}--\ref{fig:upperistight} show $\gamma_N^\star$ along with upper bounds, $\overline \gamma_N$, and lower bounds, $\underline \gamma_N$ for four different pairs of scalar systems, defined in Table~\ref{tab:parameters}.
The optimal performance level, $\gamma_N^\star$, was computed using the construction in Appendix~\ref{app:exact}, gridding the probability simplex $\{(\theta, 1-\theta) : \theta = 0, 10^{-3}, \ldots, 1-10^{-3}, 1\}$ and the bounds were computed using Theorem~\ref{thm:suff} and~\ref{thm:nec}, bisecting over $\gamma$ to an accuracy of $\pm 10^{-3}$.
The systems in Fig.~\ref{fig:unstablepm} are unstable and \emph{indistinguishable}, and the resulting optimal performance level $\gamma_N^\star$ grows exponentially in $N$.
Fig.~\ref{fig:stablepm} is also \emph{indistinguishable}, but here both systems are stable. 
The optimal performance level $\gamma_N^\star$ is bounded and is equal to the lower bound $\underline \gamma_N$.
This is because the systems are BIBO stable, so picking $\hat x_N = 0$ results in an estimation error bounded by the disturbance's norm.
Fig.~\ref{fig:nice} contains two stable systems that are \emph{distinguishable}.
The performance level $\gamma^\star$ is similar to the case where the system is known, and the bounds are close. $\gamma^\star_N$ is smaller than the other examples.
Fig.~\ref{fig:upperistight} contains two unstable \emph{distinguishable} systems. Here $\gamma_N^\star$ is bounded and approaches the upper bound $\overline \gamma_N$. 
\begin{table}
	\centering
	\caption{Parameters for the systems in Fig.~\ref{fig:unstablepm}--\ref{fig:upperistight}. 
		In all cases $Q_1 = Q_2 = R_1 = R_2 = 1$ and $P_{0, i}$ is the stationary solution to~\eqref{eq:Kalman_gain}.
	} \label{tab:parameters}
	\begin{tabular}{lrrrrrr}
		\toprule
		\textbf{System} & $F_1$ & $F_2$ & $H_1$ & $H_2$ & $P_1$ & $P_2$ \\
		\midrule
		\ref{fig:unstablepm} & 1.1 & 1.1 & 1 & -1 & 1.77 & 1.77 \\
		\ref{fig:stablepm} & 0.9 & 0.9 & 1 & -1 & 1.48 & 1.48 \\
		\ref{fig:nice} & 0.7 & 0.9 & 1.5 & 1 & 1.16 & 1.48 \\
		\ref{fig:upperistight} & 2 & 1 & 1 & 16 & 4.23 & 1.00 \\
		\bottomrule
	\end{tabular}
\end{table}
\forECC{
\begin{figure*}
    \centering 
  \begin{subfigure}{0.24\linewidth}
    \centering
    \resizebox{\linewidth}{!}{\input{\figpath/fig1}}
    \caption{Two \emph{unstable} indistinguishable systems.} \label{fig:unstablepm}
  \end{subfigure}%
  \hfill
  \begin{subfigure}{0.24\linewidth}
    \centering
    \resizebox{\linewidth}{!}{\input{\figpath/fig2}}
    \caption{Two \emph{stable} indistinguishable systems.}\label{fig:stablepm}
\end{subfigure}%
  \hfill
  \begin{subfigure}{0.24\linewidth}
    \centering
    \resizebox{\linewidth}{!}{\input{\figpath/fig3}}
    \caption{Two \emph{stable} distinguishable systems.} \label{fig:nice}
  \end{subfigure}
  \hfill
  \begin{subfigure}{0.24\linewidth}
    \centering
    \resizebox{\linewidth}{!}{\input{\figpath/fig4}}
    \caption{Two \emph{unstable} distinguishable systems.} \label{fig:upperistight}
  \end{subfigure} 
  \caption{Numerically evaluated optimal performance levels, upper and lower bounds for the four system pairs considered in Section~\ref{sec:examples}. Only stable \emph{and or} distinguishable systems have bounded performance levels. In two pairs $\gamma_N^\star$ achieves the lower bound, and in Fig~\ref{fig:upperistight} it approaches the upper bound.}
\end{figure*}
}
\forArxiv{
\begin{figure*}
    \centering 
  \begin{subfigure}{0.48\linewidth}
    \centering
    \resizebox{\linewidth}{!}{\input{\figpath/fig1}}
    \caption{Two \emph{unstable} indistinguishable systems.} \label{fig:unstablepm}
  \end{subfigure}%
  \hfill
  \begin{subfigure}{0.48\linewidth}
    \centering
    \resizebox{\linewidth}{!}{\input{\figpath/fig2}}
    \caption{Two \emph{stable} indistinguishable systems. \\}\label{fig:stablepm}
\end{subfigure}%

  \begin{subfigure}{0.48\linewidth}
    \centering
    \resizebox{\linewidth}{!}{\input{\figpath/fig3}}
    \caption{Two \emph{stable} distinguishable systems.} \label{fig:nice}
  \end{subfigure}
  \hfill
  \begin{subfigure}{0.48\linewidth}
    \centering
    \resizebox{\linewidth}{!}{\input{\figpath/fig4}}
    \caption{Two \emph{unstable} distinguishable systems.} \label{fig:upperistight}
  \end{subfigure} 
  \caption{Numerically evaluated optimal performance levels, upper and lower bounds for the four system pairs considered in Section~\ref{sec:examples}. Only stable \emph{and or} distinguishable systems have bounded performance levels. In two pairs $\gamma_N^\star$ achieves the lower bound, and in Fig~\ref{fig:upperistight} it approaches the upper bound.}
\end{figure*}
}
\section{Conclusions} \label{sec:conclusions}
This article proposed a method to compute upper and lower bounds for the optimal minimax performance level for uncertain linear systems, where the uncertainty belongs to a finite set.
The bounds are computed by evaluating the positive-definiteness \onote{of matrices} appearing in coupled Riccati recursions.
The performance level refines the notion of \emph{distinguishability} in a priori analysis of the problem set-up for multiple-model estimation, \onote{and answers} the question ``To what extent can I guarantee the performance multiple-model estimation applied to my problem?''.
Our experiments indicate that if similar output trajectories come from similar state trajectories, the gain is small.
This agrees with the intution that such systems generate similar estimates, and that in order for these estimates to be poor, the disturbances must be large.
However, if similar output trajectories come from different state trajectories, the state estimates will be different even for small disturbances, and as the optimal estimate is an interpolation of the estimates from the different models, the term $x - \hat x_N$ will be large even for small distrubances.
\onote{The provided examples show that there are systems where the optimal performance level is equal to its lower bound, approaches its upper bound, and where neither bound is ever tight.}

As with $\mathcal H_\infty$-control and estimation, the results are valid for any disturbance realization but are conservative if good disturbance statistics are available. 

\subsection{Future Work}
The numerical examples show that the bounds are tight for some systems, but not for others. 
The difference between the upper and lower bounds trivially bounds the conservativeness, but obtaining general conditions, and classifying systems where the bounds are tight, would enhance the practical utility of the results.

In this work, the system parameters $F_i$ and $H_i$ are assumed to be fixed. 
The extension to time-varying parameters is straightforward, but the extension to jump-linear systems is not.
The reason is that the number of feasible parameter trajectories grows exponentially with time.
There are heuristic ways of combining the Kalman filter estimates from different models, such as Blom's interacting-multiple-model estimator,~\cite{Blom1984}.

The worst-case history can be losslessly compressed to quadratic functions, but the number of functions will grow exponentially in time.
However, it is possible to upper bound the time-evolution of the worst-case \onote{data-consistent} parameter realization by updating a constant number of quadratic functions, similar to how we combine many Kalman-filter estimates into one estimate in this paper.
It would be interesting to exploit this bound to extend the results to jump-linear systems.

%% file: figures/fig1.tex
\pgfplotstableread[col sep=comma]{\figpath/fig1.csv}{\datatable}
\begin{tikzpicture}
    \begin{axis}[%
        xlabel={$N$},
        ylabel={$\gamma_N$},
        width={5cm},
        height={5cm},
        axis x line={bottom},
        axis y line={left},
        xminorticks = {false},
        yminorticks = {false},
	legend pos=north west
    ]
	\addplot[line width = 1pt, blue] table[x=t, y=gammas] from \datatable;
	\addlegendentry{$\gamma^\star_N$}
	\addplot[line width = 2pt, loosely dotted, black] table[x=t, y=lowers] from \datatable;
	\addlegendentry{$\underline\gamma_N$}
	\addplot[line width = 2pt, dash pattern = on 8pt off 4pt, black] table[x=t, y=uppers] from \datatable;
	\addlegendentry{$\overline\gamma_N$}
    \end{axis}
\end{tikzpicture}

%% file: figures/fig2.tex
\pgfplotstableread[col sep=comma]{\figpath/fig2.csv}{\datatable}
\begin{tikzpicture}
    \begin{axis}[%
        xlabel={$N$},
        ylabel={$\gamma_N$},
        width={5cm},
        height={5cm},
        axis x line={bottom},
        axis y line={left},
        xminorticks = {false},
        yminorticks = {false},
    ]
	\addplot[line width = 1pt, blue] table[x=t, y=gammas] from \datatable;
	\addplot[line width = 2pt, loosely dotted, black] table[x=t, y=lowers] from \datatable;
	\addplot[line width = 2pt, dash pattern = on 8pt off 4pt, black] table[x=t, y=uppers] from \datatable;
    \end{axis}
\end{tikzpicture}

%% file: figures/fig3.tex
\pgfplotstableread[col sep=comma]{\figpath/fig3.csv}{\datatable}
\begin{tikzpicture}
    \begin{axis}[%
        xlabel={$N$},
        ylabel={$\gamma_N$},
        width={5cm},
        height={5cm},
        axis x line={bottom},
        axis y line={left},
        xminorticks = {false},
        yminorticks = {false},
	ymin=1.21,
	ymax=1.36,
    ]
	\addplot[line width = 1pt, blue] table[x=t, y=gammas] from \datatable;
	\addplot[line width = 2pt, loosely dotted, black] table[x=t, y=lowers] from \datatable;
	\addplot[line width = 2pt, dash pattern = on 8pt off 4pt, black] table[x=t, y=uppers] from \datatable;
    \end{axis}
\end{tikzpicture}

%% file: figures/fig4.tex
\pgfplotstableread[col sep=comma]{\figpath/fig4.csv}{\datatable}
\begin{tikzpicture}
    \begin{axis}[%
        xlabel={$N$},
        ylabel={$\gamma_N$},
        width={5cm},
        height={5cm},
        axis x line={bottom},
        axis y line={left},
        xminorticks = {false},
        yminorticks = {false},
	xmax=20,
    ]
	\addplot[line width = 1.6pt, blue] table[x=t, y=gammas] from \datatable;
	\addplot[line width = 2pt, loosely dotted, black] table[x=t, y=lowers] from \datatable;
	\addplot[line width = 2pt, dash pattern = on 8pt off 4pt, black] table[x=t, y=uppers] from \datatable;
    \end{axis}
\end{tikzpicture}

%% file: tex/appendix.tex
\appendix
\forArxiv{
	\section*{Appendix}
}
\forArxiv{\section{Proofs}}
\forECC{\subsection{Proofs}}
\label{sec:proofs}
This section proves Theorems \ref{thm:suff} and \ref{thm:nec}.
In doing so we obtain an expression that can be used to evaluate the value~\eqref{eq:def_J}, but is computationally intractable for problems with uncertainties belonging to moderately-sized sets.

\forArxiv{\subsection{Proof strategy}}
\forECC{\subsubsection{Proof strategy}}
We reparameterize the disturbance trajectory $(\seq{w}{0}{N-1}, \seq{v}{0}{N-1})$ in the state-output trajectory and the active model $(\seq{x}{0}{N-1}, \seq{y}{0}{N-1}, i)$.
This reparameterization allows us to partially switch the order of the minimization and the maximization, as $\mu$ is a function of $\seq{y}{0}{N-1}$, yielding a problem of the form $\max_{\seq{y}{0}{N-1}}\min_\mu \max_{i, \seq{x}{0}{N}}$.
Previous work, \cite{Kjellqvist2022Minimax}, shows how to maximize over $\seq{x}{0}{N}$ using forward dynamic programming, resulting in the forward Riccati recursions~\eqref{eq:Kalman_gain}.

We then reformulate the maximization over the feasible set to maximizing over its convex hull.
This reformulation allows us to switch the order of minimizing with respect to $\mu$ and maximizing with respect to the model.
The catch is that while the value is unchanged, the maximizing $\theta$ is not necessarily the same.
As we are interested in the value, we can ignore this issue.

The inner minimization problem is unconstrained and convex-quadratic in the estimate $\hat x_N$, which has a closed-form solution.
The maximization over the convex hull of the model set is then bounded from above and from below by a maximum over a finite number of functions that linear-quadratic regulator costs in $\seq{y}{0}{N-1}$, which has a solution expressed by the backward Riccati recursion, \eqref{eq:riccati}.

\forArxiv{\subsection{Reparameterization}}
\forECC{\subsubsection{Reparameterization}}
The disturbance $w_t$ is uniquely determined by $F = F_i$ and $(x_{t+1}, x_t)$, and $v_t$ is uniquelly determined by $H = H_i$, $y_t$ and $x_t$.
As the maximizing player is aware of the dynamics, $i$, we can substitute $w_t = x_{t+1} - F_ix_t$ and $v_t = y_t - H_ix_t$ into~\eqref{eq:def_J},
\begin{multline}\label{eq:J_y}
    J_N^\star(\hat x_0) = \inf_\mu\sup_{\seq{x}{0}{N}, \seq{y}{0}{N-1}, i}\Bigg\{
        |x_N - \hat x_N|^2 - \gamma^2|x_0 - \hat x_{0, i}|^2_{P_{0, i}}  \\
         - \gamma^2\sum_{t=0}^{N-1}\Big[|x_{t + 1} - F_ix_t|^2_{Q_i^{-1}} + |y_t - H_ix_t|^2_{R_i^{-1}} \Big]        \Bigg\}.
\end{multline}
Furthermore, as $\mu$ is a function of $\seq{y}{0}{N-1}$, we can move the maximization over output trajectories outside of the minimization and minimize directly over the estimate $\hat x_N \in \R^n$.
Consider the inner maximization over state trajectories, which is a function of the observations and estimates,
\begin{multline} \label{eq:inner}
	J^\text{inner}_N(\seq{y}{0}{N-1}, \hat x_N, \hat x_0) = \forECC{\\}
    \sup_{\seq{x}{0}{N}, i} \Bigg\{
        |x_N - \hat x_N|^2  - \gamma_N^2|x_0 - \hat x_{0, i}|^2_{P_{0, i}}\\
         - \gamma^2\sum_{t=0}^{N-1}\Big[|x_{t + 1} - F_ix_t|^2_{Q_i^{-1}} + |y_t - H_ix_t|^2_{R_i^{-1}} \Big]        \Bigg\}.
\end{multline}
Then~\eqref{eq:J_y} can be written as
\begin{equation} \label{eq:J_y_mu_x}
    J_N^\star(\hat x_0) = \sup_{\seq{y}{0}{N-1}}\inf_{\hat x_N}J^\text{inner}_N(\seq{y}{0}{N-1}, \hat x_N, \hat x_0).
\end{equation}

\forArxiv{\subsection{Foward recursion}}
\forECC{\subsubsection{Foward recursion}}
Following the proof of Theorem 1 in \cite{Kjellqvist2022Minimax}, with $P_{0, i}$ as the stationary solution to~\eqref{eq:Kalman_gain}, we see that the value inner optimization problem \eqref{eq:inner} is equal to
\begin{align}
	& \sup_{i, x_N}\left\{|\hat x_N - x_N|^2 - \gamma^2 \left(|x_N - \breve x_{N, i}|^2_{P^{-1}_i} + c_{N, i} \right) \right\} \nonumber \\
	& \quad = \max_i \left\{|\hat x_N - \breve x_{N, i}|^2_{(I - \gamma^{-2}P_i)^{-1}} - \gamma^2 c_{N, i} \right\}, \label{eq:equivalent_inner}
\end{align}
if $I \succ \gamma^{-2} P_i^{-1}$ for all $i$. The value is unbounded if $I \nsucceq \gamma^{-2} P_i^{-1}$ for some $i$, which proves Proposition~\ref{prop:bound_P}.
Proposition~\ref{prop:minimax_optimal} shows how to compute $\breve x_{N, i}$, $P_i$ and $c_{N, i}$ in \eqref{eq:equivalent_inner}

\subsection{Exact computations of $J^\star_N$}\label{app:exact}
By substituting \eqref{eq:equivalent_inner}, we see that the value of \eqref{eq:J_y_mu_x} is equal to
\begin{equation}\label{eq:equivalent_cost}
    \sup_{\seq{y}{0}{N-1}}\min_{\hat x_N}\max_i \Big\{|\hat x_N - \breve x_{N, i}|^2_{(I - \gamma^{-2}P_i)^{-1}} - \gamma^2 c_{N, i} \Big\}.
\end{equation}

Maximizing over the finite set $\modelset$ in \eqref{eq:equivalent_cost} is equivalent to optimizing for convex combinations over the probability simplex $\Theta = \{\theta \in \R^n :0 \leq  \theta_i \leq 1, \sum_{i=1}^\nmodels \theta_i = 1\}$.
The equivalence is because the optimal value of a linear program over a simplex is located on a vertix. 
As \eqref{eq:equivalent_cost} is convex in $\hat x$, the minimizing $\hat x$ can be bounded in terms of $\breve x_{N, i}$.
The convex combination is affine in $\theta$, so Von Neumann's minmax theorem applies and the value \eqref{eq:equivalent_cost} is equal to
\begin{equation*}
	\sup_{\theta \in \Theta}\min_{\hat x}\Bigg\{
		\sum_{i=1}^\nmodels\theta_i \left( |\hat x - \breve x_{N, i}|^2_{Q_{N, i}} - \gamma^2c_{N, i}\right)
	\Bigg\},
\end{equation*}
where $Q_{N, i} = (I - \gamma^{-2}P_i)^{-1}$.
Applying Lemma~\ref{lemma:interpolation} to the inner minimization problem means that the value \eqref{eq:equivalent_cost} is equal to
\forECC{%
\begin{multline}\label{eq:Jytheta}
 \sup_{\seq{y}{0}{N-1}, \theta} \Bigg\{
 \sum_{i = 1}^\nmodels\theta_i\left( |\breve x_{i, i}|^2_{Q_{N, i}} - \gamma^2c_{N, i}\right) \\
		-\Big|\sum_{i = 1}^\nmodels\theta_i Q_{N, i}\breve x_{N, i}\Big|^2_{(\sum \theta_i Q_{N, i})^{-1}}
	\Bigg\}.
\end{multline}
}
\forArxiv{%
\begin{equation}\label{eq:Jytheta}
 \sup_{\seq{y}{0}{N-1}, \theta} \Bigg\{
 \sum_{i = 1}^\nmodels\theta_i\left( |\breve x_{i, i}|^2_{Q_{N, i}} - \gamma^2c_{N, i}\right) 
		-\Big|\sum_{i = 1}^\nmodels\theta_i Q_{N, i}\breve x_{N, i}\Big|^2_{(\sum \theta_i Q_{N, i})^{-1}}
	\Bigg\}.
\end{equation}
}
For a fixed $\theta$, this is a sequential quadratic optimization problem in $y$ that can be solved using dynamic programming.
In fact this can be reformulated into a standard linear-quadratic regulator problem, except that the terminal penalty is indefinite.
This indefinite term will, for small values of $\gamma_N$, lead to a loss of concavity in $\seq{y}{0}{N-1}$.
This means that the value is unbounded, and $\gamma_N < \gamma_N^\star$.
Larger values of $\gamma_N$ will compensate for the indefinite term and ensure concavity in $\seq{y}{0}{N-1}$.
Testing for concavity amounts to evaluating whether $\stack X_t$ in \eqref{eq:Riccati_bw} is positive definite for all $t$.
If concavity in $\seq{y}{0}{N-1}$ holds for all $\theta \in \Theta$, then the value is finite and $\gamma_N \geq \gamma_N^\star$.
Define
\begin{equation*}
\begin{split}
	&\stack F \defeq \bdiag\left( \{F_i - K_i H_i\}_{i=1}^\nmodels\right) \\
	&\stack{\breve x_t} \defeq \bmat{\breve x_{t, 1}^\tran & \cdots & \breve x_{t, \nmodels}^\tran}^\tran, \quad \stack K \defeq \bmat{K_{1}^\tran & \cdots & K_{\nmodels}^\tran}^\tran.
\end{split}
\end{equation*}
Then, the multi-observer update \eqref{eq:estimators:state_update} becomes,
\[
	\stack{\breve x}_{t + 1} = \stack F\stack{\breve x_t} + \stack Ky_t.
\]
Further, let 
\begin{align}
	& \stack Q_N  \defeq \bdiag\{\theta_i Q_{N, i}\}_{i=1}^\nmodels \forECC{\nonumber \\ 
	& \quad} -
	\bmat{
		\theta_{1}Q_{N, 1}\\
		\vdots \\
		\theta_{\nmodels}Q_{N, \nmodels}
	}
	\left( (\sum \theta_i Q_{N, i})^{-1}\right)^{-1}
	\bmat{
		\theta_{1}Q_{N, 1}\\
		\vdots \\
		\theta_{\nmodels}Q_{N, \nmodels}
	}^\tran, \\ \nonumber
		   & \bmat{
		\stack Q & \stack N^\tran \\
		\stack N & \stack R
	}
		   \defeq 
	\left[
		\begin{array}{c}
			\bdiag\left(\{-H_i^\tran\}_{i=1}^\nmodels\right) \\
			\begin{array}{ccc}
				I & \cdots & I
			\end{array}
		\end{array}
\right]\nonumber \\
		   & \quad \times\bdiag\left(\{\theta_i \tilde R_i^{-1}\}_{i=1}^\nmodels\right)
	\left[
		\begin{array}{c}
			\bdiag\left(\{-H_i^\tran\}_{i=1}^\nmodels\right) \\
			\begin{array}{ccc}
				I & \cdots & I
			\end{array}
		\end{array}
\right]^\tran,\label{eq:one_step:weights}
\end{align}
where $\times$ denotes standard matrix product.
With 
\forECC{$l(\theta, \stack{\breve x}_t, y_t) = \gamma^2_N \left(|\stack{\breve x}_t|^2_{\stack Q} - 2 y_t^\tran \stack N\stack{\breve x_t} + |y_t|^2_{\stack R}\right)$,}
\forArxiv{$$l(\theta, \stack{\breve x}_t, y_t) = \gamma^2_N \left(|\stack{\breve x}_t|^2_{\stack Q} - 2 y_t^\tran \stack N\stack{\breve x_t} + |y_t|^2_{\stack R}\right),$$}
\eqref{eq:J_y_mu_x} becomes
\begin{equation}\label{eq:one_step:lqr}
J^\star_N(\hat x_0) = -\inf_{\theta}\underbrace{\inf_{\seq{y}{0}{N-1}}\left\{ 
		|\stack {\breve x_t}|^2_{\stack Q_N} + \sum_{t=0}^{N-1} l(\theta, \stack{\breve x}_t, y_t)
\right\}}_{\defeq J_N(\theta, \hat x_0)}.
\end{equation}

It is apparent that $l$ is strictly convex in $y_t$.
However, the terminal penalty matrix, $\stack Q_N$, is indefinite, which may cause \eqref{eq:one_step:lqr} to lose convexity and become unbounded.
\begin{remark}
	The stage cost is a convex combination of the Kalman filter residuals $l(\theta, \stack{\breve x}_t, y_t) = \gamma_N^2\sum_{i=1}^\nmodels\left( \theta_i c_{t, i}\right)$.
\end{remark}
The Riccati recursions corresponding to the linear-quadratic regulator are well described in many textbooks, for instance in~\cite[Chapter 11.2]{Astrom97}, and can be used to compute the value provided that $\theta \in \Theta$ fixed:
\begin{equation}\label{eq:Riccati_bw}
	\begin{aligned}
		\stack X_t & = \stack K^\tran \stack T_t \stack K + \stack R, \quad
		\stack L_t = \stack X_t^{-1} (\stack K^\tran \stack T_t \stack F - \stack N) \\
		\stack T_{t-1} & = \stack F^\tran \stack T_t \stack F + \stack Q - \stack L_t^\tran \stack X_t \stack L_t. \\
	\end{aligned}
\end{equation}
The relationship between the solution to the above Riccati equations and the value of the game are summarized in the below lemma.
\begin{lemma}\label{lemma:lqr}
	Consider the backward Riccati equations above with terminal condition $\stack T_N = -\stack Q_N/ \gamma_N^2$.
	Let $J_N^\star(x_0)$ be the value of the game \eqref{eq:def_J} and $J_N(\theta, x_0)$ be value of the inner, sequantial, optimization problem in~\eqref{eq:one_step:lqr}.
	If $\stack X_t \not \preceq 0$ for some $\theta \in \Theta$, then $J_N^\star(x_0)$ is unbounded.
	If $\stack X_t \succ 0$ for all $\theta \in \Theta$ then $J_N(\theta, \hat x_0) = -\gamma_N^2 |\stack{\breve x}_0|^2_{\stack T_0}$, and
	\[
		J_N^\star(\hat x_0) = \max_\theta \left(J_N(\theta, \hat x_0) \right).
	\]
\end{lemma}

\forArxiv{\section{Upper- and Lower bounds of $J^\star_N$}}
\forECC{\subsection{Upper- and Lower bounds of $J^\star_N$}}
This section develops  upper and lower bounds on the objective, \eqref{eq:def_J}.
As the maximum is greater than the average of any two points, we have that
\begin{align}
    J_N^\star(\hat x_0) & \geq \sup_{i, j, \seq{y}{0}{N-1}}\min_{\hat x_N} \frac{1}{2}\Big\{
        |\hat x_N - \breve x_{N, i}|^2_{(I - \gamma^{-2}P_i)^{-1}} \nonumber \\
        & \quad - \gamma^2 c_{N, i} + |\hat x_N - \breve x_{N, j}|^2_{(I - \gamma^{-2}P_j)^{-1}} - \gamma^2 c_{N, j} 
     \Big\} \nonumber \\
        & =\sup_{i, j, \seq{y}{0}{N-1}}  \frac{1}{2}\Big \{
            |\breve x_{N, i} - \breve x_{N, j}|^2_{(2I - \gamma^{-2}P_i - \gamma^{-2}P_j)^{-1}} \forECC{\nonumber \\
	& \quad} - \gamma^2c_{N, i} - \gamma^2c_{N, j}
    \Big \} \forArxiv{\nonumber \\ &}\defeq \max_{i, j}\underline{J}^{ij}_N(\hat x_0). \label{eq:lower_bound}
\end{align}
Thus $\gamma_N < \gamma_N^\star$ only if $\underline J^{ij}_N(\hat x_0)$ is bounded for all pairs $(i, j)$.
Towards finding a sufficient condition, let $S\in\R^{n \times n}$ be a positive definite matrix such that $ S \preceq I - \gamma^{-2}P_i$ for all $i = 1, \ldots, \nmodels$. 
Then, applying Lemma~\ref{lemma:quad_prog} to \eqref{eq:Jytheta}, we have
\begin{align}
& J_N^\star(\hat x_0)  \leq  \forECC{\nonumber \\
& \quad} \max_{y,\theta}\Big\{ \sum_{i, j}^{\nmodels}\theta_i \theta_j |\breve x_{N, i} - \breve x_{N, j}|^2_{S^{-1}}/2- \gamma_N^2\sum_i \theta_i c_{N, i} \Big\} \nonumber\\
	\forArxiv{ & \quad \leq \frac{1}{2}\max_y \max_{\theta}\sum_i^\nmodels\theta_i\Bigg[- \gamma^2 c_{N, i}\forECC{\nonumber\\
		   & \qquad}  + \max_\sigma\Big\{\sum_j^\nmodels\sigma_j(|\breve x_{N, i} - \breve x_{N, j}|^2_{S^{-1}}- \gamma^2 c_{N, j})\Big\} \Bigg]\nonumber \\}
		   & \quad \forArxiv{=} \forECC{\leq} \max_{i, j}\underbrace{\max_y\frac{1}{2} \left \{|\breve x_{N, i} - \breve x_{N, j}|^2_{S^{-1}} - \gamma^2 (c_{N, i} + c_{N, j})  \right\}}_{\overline{J}^{ij}_N(\hat x_0)}. \label{eq:upper_bound}
\end{align}
Thus, if $\overline{J}^{ij}_N(\hat x_0)$ is bounded for all pairs $(i, j)$, then $\gamma^\star_N \leq \gamma_N$.
The only difference between the expressions of $\overline J^{ij}(\hat x_0)$ and $\underline J^{ij}(\hat x_0)$ is the penalty of the term $|\breve x_{N, i} - \breve x_{N, j}|^2_*$.

Theorems~\ref{thm:suff} and \ref{thm:nec} follow from applying Lemma~\ref{lemma:lqr} to the upper bound $\overline J^{ij}_N(\hat x_0)$ in~\eqref{eq:upper_bound} and the lower bound $\underline J^{ij}_N(\hat x_0)$ in~\eqref{eq:lower_bound} with $\theta_i = \theta_j = \frac{1}{2}$.

\forArxiv{\section{Lemmata}}
\forECC{\subsection{Lemmata}}
\begin{lemma} \label{lemma:quad_prog}
    Let $X_i \succ 0$ and $\theta_i \in (0, 1)$ for $i = 1, \ldots, \nmodels$ and that $\sum_{i=1}^\nmodels \theta_i = 1$.
    Let $S = \sum \theta_{i=1}^\nmodels X_i^{-1}$, then
    \begin{multline*}
        \min_v\left\{ \sum \theta_i|v - x_i|^2_{X_i^{-1}} \right\}
	= \sum_{i=1}^\nmodels\theta_i\Big(|X_i^{-1}x_i|^2_{X_i - S^{-1}} \\
	+ \frac{1}{2}\sum_{j=1}^\nmodels\theta_j \left(|X_i^{-1} x_i - X_j^{-1}x_j|^2_{S^{-1}}\right) \Big).
    \end{multline*}
\end{lemma}
\begin{proof}
    As $X_i \succ 0$, we have that $\sum \theta_i X_i^{-1} \succ 0$ and the (unique) minimum is a stationary point.
    We have
    \forECC{%
    \begin{multline*}
        \min_v\left\{ \sum \theta_i|v - x_i|^2_{X_i^{-1}} \right\} = \sum_{i = 1}^\nmodels \theta_i|x_i|^2_{X_i^{-1}}\\
	- |\sum_{i=1}^\nmodels \theta_iX_i^{-1}x_i|^2_{(\sum_1^\nmodels\theta_iX_i^{-1})^{-1}}
\end{multline*}
} \forArxiv{%
    \begin{equation*}
        \min_v\left\{ \sum \theta_i|v - x_i|^2_{X_i^{-1}} \right\} = \sum_{i = 1}^\nmodels \theta_i|x_i|^2_{X_i^{-1}}
	- |\sum_{i=1}^\nmodels \theta_iX_i^{-1}x_i|^2_{(\sum_1^\nmodels\theta_iX_i^{-1})^{-1}}
\end{equation*}
}

    With $S := (\sum_1^\nmodels \theta_i X_i^{-1})$, we have that
    \begin{align*}
       & -|\sum_{i=1}^\nmodels \theta_iX_i^{-1}x_i|^2_S = -\sum_{i=1}^\nmodels\sum_{j=1}^\nmodels \theta_i \theta_j x_i^\top X_i^{-\top}SX_j^{-1}x_j \\
                        &  \quad=\frac{1}{2}\sum_{i=1}^\nmodels\sum_{j=1}^\nmodels\theta_i \theta_j \left(|X_i^{-1} x_i - X_j^{-1}x_j|^2_S\right) - \sum_{i=1}^\nmodels \theta_i |X_i^{-1}x_i|^2_S.
    \end{align*}
\end{proof}

\begin{lemma}[Interpolation]\label{lemma:interpolation}
	Let $z_k \in \R^n$ and $Z_k \in \R^{n \times n}$ be matrices such that $\sum_{k = 1}^{K}Z_k \succ 0$ for $k = 1, \ldots, K$. 
	Then,
	\[
		\min_x\left\{\sum_{k = 1}^K|x - z_k|^2_{Z_k} \right\} = \sum_{k = 0}^K|z_k|^2_{Z_k} - \left |\sum_{k=1}^K Z_k z_k\right|^2_{(\sum_{k = 1}^K Z_k)^{-1}}.
	\]
\end{lemma}
\begin{proof}
	The problem is unconstrained and strictly convex---the minimizing solution is given by the stationary point.
\end{proof}